\documentclass[11pt]{amsart}
\usepackage{latexsym}
\usepackage{amscd}

\normalsize

\RequirePackage{xspace}

\newcommand{\angles}[1]         {{\langle #1 \rangle}}

\newcommand{\thought}[1]{}
\renewcommand{\thought}[1]{ \textbf{[#1]}}

\usepackage{enumerate}        
\newenvironment{roenumerate}{\begin{enumerate}[\upshape (i)]}{\end{enumerate}}

\newcommand{\iref}[1]{(\ref{#1})}

\newcommand{\df}[1]{\emph{#1}}

\newcommand{\period}    {{\makebox[0pt][l]{\hspace{2pt} .}}}

\newcommand\nc {\newcommand}
\newcommand\rnc{\renewcommand}
\font\tenfrakt=eufm10   \font\tenscript=eusm10
\font\sevenfrakt=eufm7  \font\sevenscript=eusm7
\font\fivefrakt=eufm5   \font\fivescript=eusm5

\newfam\fraktfam
\textfont\fraktfam=\tenfrakt
\scriptfont\fraktfam=\sevenfrakt
\scriptscriptfont\fraktfam=\fivefrakt

\newfam\scriptfam
\textfont\scriptfam=\tenscript
\scriptfont\scriptfam=\sevenscript
\scriptscriptfont\scriptfam=\fivescript
\def\script{\fam\scriptfam\tenscript}
\newcount\blopone
\newcount\xone
\newcount\xtwo
\newcount\ytwo

\newtheorem{theorem}{Theorem}[section]
\newtheorem{prop}[theorem]{Proposition}

\newtheorem{refinement}[theorem]{Refinement}
\newtheorem{summary}[theorem]{Summary}
\newtheorem{importnota}[theorem]{Important Notation}
\newtheorem{prblm}[theorem]{Problem}
\newtheorem{notation}[theorem]{Notation}
\newtheorem{defin}[theorem]{Definition}
\newtheorem{caution}[theorem]{Caution}
\newtheorem{remark}[theorem]{Remark}
\newtheorem{reminder}[theorem]{Reminder}
\newtheorem{lemma}[theorem]{Lemma}
\newtheorem{construction}[theorem]{Construction}
\newtheorem{corollary}[theorem]{Corollary}
\newtheorem{example}[theorem]{Example}
\newtheorem{conclusion}[theorem]{Conclusion}
\newtheorem{triviality}[theorem]{Triviality}
\newtheorem{proto}[theorem]{Prototype Quasifibration}
\newtheorem{cauex}[theorem]{Cautionary Example}
\newtheorem{hypo}[theorem]{Hypothesis}
\newtheorem{subth}{ }[theorem]
\newtheorem{case}{Case}[theorem]
\newtheorem{ssubth}{ }[subth]

\nc\tri[1]{\begin{triviality}
\label{#1}}
\nc\cas[1]{\begin{case}
\label{#1}
\begin{em}}
\nc\rfn[1]{\begin{refinement}
\label{#1}}
\nc\prt[1]{\begin{proto}
\label{#1}}
\nc\lem[1]{\begin{lemma}
\label{#1}}
\nc\pro[1]{\begin{prop}
\label{#1}}
\nc\thm[1]{\begin{theorem}
\label{#1}}
\nc\cor[1]{\begin{corollary}
\label{#1}}
\nc\dfn[1]{\begin{defin}
\label{#1}}
\nc\sthm[1]{\begin{subth}
\label{#1}}
\nc\exm[1]{\begin{example}
\label{#1}
\begin{em}}
\nc\plm[1]{\begin{prblm}
\label{#1}
\begin{em}}
\nc\rmk[1]{\begin{remark}
\label{#1}
\begin{em}}
\nc\rmd[1]{\begin{reminder}
\label{#1}
\begin{em}}
\nc\ntn[1]{\begin{notation}
\label{#1}
\begin{em}}
\nc\smr[1]{\begin{summary}
\label{#1}
\begin{em}}
\nc\cau[1]{\begin{caution}
\label{#1}
\begin{em}}
\nc\hyp[1]{\begin{hypo}
\label{#1}}
\nc\imn[1]{\begin{importnota}
\label{#1}
\begin{em}}
\nc\cax[1]{\begin{cauex}
\label{#1}
\begin{em}}
\nc\con[1]{\begin{construction}
\label{#1}
\begin{em}}
\nc\ssthm[1]{\begin{ssubth}
\label{#1}
\begin{em}}
\nc\cnc[1]{\begin{conclusion}
\label{#1}
\begin{em}}

\nc\elem{\end{lemma}}
\nc\erfn{\end{refinement}}
\nc\eprt{\end{proto}}
\nc\ethm{\end{theorem}}
\nc\ecor{\end{corollary}}
\nc\edfn{\end{defin}}
\nc\esthm{\end{subth}}
\nc\epro{\end{prop}}
\nc\etri{\end{triviality}}
\nc\eexm{\end{em}
\end{example}}
\nc\ermk{\end{em}
\end{remark}}
\nc\ermd{\end{em}
\end{reminder}}
\nc\eplm{\end{em}
\end{prblm}}
\nc\ecas{\end{em}
\end{case}}
\nc\ecau{\end{em}
\end{caution}}
\nc\ecax{\end{em}
\end{cauex}}
\nc\eimn{\end{em}
\end{importnota}}
\nc\entn{\end{em}
\end{notation}}
\nc\econ{\end{em}
\end{construction}}
\nc\esmr{\end{em}
\end{summary}}
\nc\ehyp{
\end{hypo}}
\nc\ecnc{\end{em}
\end{conclusion}}
\nc\essthm{\end{em}
\end{ssubth}}

\nc\sst{\scriptstyle}
\newcommand{\comment}[1]{}
\newcommand{\ri}{\longrightarrow}
\newcommand{\sr}{\rightarrow}

\newcommand{\zz}{{\mathbb Z}}

\newcommand{\C}{{\mathbb C}}

\nc\z{\zeta}
\nc\bc{{\mathbb{BC}}}
\nc\ct{{\script T}}
\nc\cs{{\script S}}
\nc\car{{\script R}}
\nc\cd{{\script D}}
\nc\cc{{\script C}}
\nc\ca{{\script A}}
\nc\ci{{\script I}}
\nc\co{{\script O}}
\nc\bd{\begin{description}}
\nc\ed{\end{description}}
\nc\ctob{{\script C}at\big(\ci^{op},\ca\big)}
\nc\clim{{\ds\mathop{\rm lim}_{\ds\longleftarrow}}}
\nc\climi{\clim^{\!i}\,}
\nc\climn{\clim^{\!n}\,}
\nc\colim{{\ds\mathop{\rm colim}_{\ds\la}}}
\nc\oa{\overline{\ca}}
\nc\s{\sigma}
\nc\ta{\tau}
\nc\os{\overline\sigma}
\nc\ot{\overline\tau}
\nc\T{\Sigma}
\nc\de[1]{{\mathop{\rm deg(#1)}}}
\nc\Ad[1]{\mathop{\rm Ad}(#1)}
\nc\ad[1]{\mathop{\rm ad}(#1)}

\def\der #1 {D\left(#1\right)}
\nc\prf{\begin{proof}}
\nc\eprf{\end{proof}}
\nc\ds{\displaystyle}

\nc\cb{{\script B}}
\nc\ab{{\script A}b}

\nc\csab{{\script C}at\big(\cs^{op},\ab\big)}
\nc\ctab{{\script C}at\Big({\{\ct^\alpha\}}^{op},\ab\Big)}
\nc\csex{{\script E}x\big(\cs^{op},\ab\big)}
\nc\ctex{{\script E}x\Big({\{\ct^\alpha\}}^{op},\ab\Big)}
\nc\sub{\qquad\subset\qquad}
\nc\ctr[1]{{\left.\ct\left(-,#1\right)\right|}_{\cs}}
\nc\ctrf[2]{{\left.\ct\left(#1,#2\right)\right|}_{\cs}}
\nc\Ctr[1]{{\left.\ct\left(-,#1\right)\right|}_{\ct^\alpha}}
\nc\Ctrf[2]{{\left.\ct\left(#1,#2\right)\right|}_{\ct^\alpha}}

\nc\la{\longrightarrow}
\nc\nin{\noindent}
\nc\cad[1]{\text{card}(#1)}
\nc\eq{\quad=\quad}
\nc\BA{\begin{array}{c}}
\nc\EA{\end{array}}
\nc\barr{
\[
\begin{array}{cccccccccccccccc}
}
\nc\earr{
\end{array}
\]
}
\nc\as[1]{{\langle S\rangle}^{#1}}
\nc\sh{\hbox{\it shift}}

\nc\yy[1]{{\left.\ct\left(-,#1\right)\right|}_{\ct^c}}
\nc\vrep[2]{{\left.\ct\left(#1,#2\right)\right|}_{\ct^\alpha}}
\nc\da{\downarrow}
\nc\Hom{{\mathop{\rm Hom}}}
\nc\End{{\mathop{\rm End}}}
\nc\Ext{{\mathop{\rm Ext}}}
\nc\mr\Modtc
\nc\PExt{{\mathop{\rm PExt}}}

\nc\ra\ri
\nc\y[1]{\mathbf{y}#1}
\nc\x[1]{\mathbf{z}#1}
\nc\Mod[1]{\ensuremath{\mathop{\textup{Mod-}#1}}\xspace}
\nc\Md {\ensuremath{\mathop{\textup{Mod}}}}
\rnc\mod[1]{\ensuremath{\mathop{\textup{mod-}#1}}\xspace}
\nc\Modtc{\Mod{\ct^c}}
\nc\pgldim[1]{\mathop{\rm pgldim}\,#1}

\begin{document}

\author{J. Daniel Christensen}
\address{Department of Mathematics\\
        Johns Hopkins University\\
        3400 N. Charles Street\\
        Baltimore, MD, USA}
\email{jdc@math.jhu.edu}

\author{Bernhard Keller}
\address{Universit\'e Paris 7\\
         UFR de Math\'ematiques\\
         Institut de Math\'ematiques\\
         UMR 7586 du CNRS \\
         Case 7012\\
         2, place Jussieu\\ 
         75251 Paris Cedex 05\\
         FRANCE}
\email{keller@math.jussieu.fr}

\author{Amnon Neeman}
\address{Center for Mathematics and its Applications \\
        School of Mathematical Sciences\\
        John Dedman Building\\
        The Australian National University\\
        Canberra, ACT 0200\\
        AUSTRALIA}
\email{Amnon.Neeman@anu.edu.au}

\title[Failure of Brown representability]
{Failure of Brown representability in derived categories}

\dedicatory{dedicated to H.~Lenzing on the occasion of his sixtieth birthday}

\begin{abstract}
Let $\ct$ be a triangulated category with coproducts,
$\ct^c\subset\ct$
the full subcategory of compact objects in $\ct$. 
If $\ct$ is the homotopy category of spectra, Adams proved
the following in \cite{Adams71}:
All homological functors
${\{\ct^c\}}^{op}\ri\ab$ are the restrictions of representable
functors on $\ct$, and all natural transformations are the
restrictions of morphisms in $\ct$. 

It has been something of a mystery, to what extent this 
generalises to other triangulated categories. 
In~\cite{Neeman4}, it was proved that Adams' theorem remains true as
long as $\ct^c$ is countable, but can fail in general.
The failure exhibited was that there can be natural transformations
not arising from maps in $\ct$.

A puzzling open problem remained: Is every homological functor
the restriction of a representable functor on $\ct$?
In a recent paper, Beligiannis~\cite{Beligiannis99} made some
progress. But in this article, we settle the problem. The answer is
no. There are examples of derived categories $\ct=D(R)$
of rings, 
and homological functors ${\{\ct^c\}}^{op}\ri\ab$ which
are not restrictions of representables.
\end{abstract}

\keywords{Brown representability, derived category, purity, pure global
dimension, hereditary ring}

\maketitle

\tableofcontents

\section*{Introduction}
\label{S0}
The introduction is written for the reader who knows about 
derived categories, but is not necessarily familiar with
previous articles by the authors and their 
close friends. We begin with a sketch
of the work done in the last 10 years, generalising results
from homotopy theory to derived categories. The experts may want to 
skip this, and go directly to Notation~\ref{N00.3} 
on page~\pageref{N00.3}.
After the very general survey, will come a much more focused one.
We will give, in some detail, the history of the results on
generalising the theorem of Brown and Adams to derived categories.
Then we will explain the two open problems, which we settle
in this article. Finally, we will give the nature of our
counterexamples. 

Let $\ct$ be a triangulated category. The representable functors
$\ct(-,X)$ are all homological; that is, they take triangles to
long exact sequences. Given a triangulated
subcategory $\cs\subseteq\ct$,
we can restrict a representable functor on $\ct$ to a functor on
$\cs$. We denote the restriction by ${\ct(-,X)\big|}_\cs$.
All such functors are clearly homological.

The most interesting version of this, is where $\ct$ is a triangulated
category with coproducts, and $\cs$ is the full subcategory $\ct^{c}$ of all 
compact objects in $\ct$. 
\dfn{D0.1}
An object $c\in\ct$ is called {\em compact,} if the functor
$\ct(c,-)$ commutes with coproducts. 
\edfn 
\nin
We remind the reader that for $\ct$ the homotopy 
category of spectra, $\ct^c\subset\ct$ is the subcategory of
finite spectra. For $\ct=D(R)$, the unbounded derived category of 
right $R$--modules, $\ct^c$ turns out to be the subcategory of
perfect complexes, that is, complexes isomorphic to finite
complexes of finitely generated projective $R$--modules.
For a more detailed discussion
of examples, where $\ct$ is the unbounded derived category
of coherent sheaves on a scheme, see Sections~1 and 2
in \cite{Neeman96}.

Since the functor ${\ct(-,X)\big|}_{\ct^c}$ plays a major role in
what follows, we adopt a shorthand for it. We will write
\[
\y X={\ct(-,X)\big|}_{\ct^c}.
\]

The subject we will be studying began with a theorem of 
Adams~\cite{Adams71}.
\thm{T00.1}{\bf (Adams, 1971)}\ \ 
Let $\ct$ be the homotopy category of spectra, and $\ct^c$ the subcategory
of finite spectra. Then any homological functor
${\{\ct^c\}}^{op}\ri\ab$
is isomorphic to\/ $\y X$, for some object $X\in\ct$.
Furthermore, any natural transformation of functors
\[
\CD
\y X@>>> \y Y
\endCD
\]
is induced by some (non--unique) map $X\ri Y$.
\ethm
\rmk{R00.2}
This theorem is usually referred to as ``Brown representability''.
The reason for this is that, 10 years earlier, Brown~\cite{Brown62} 
proved a special case. In Brown's theorem, there was a
countability hypothesis on the functor. 

Calling this theorem ``Brown representability'' is somewhat confusing,
since in the same paper, Brown proved another result, somewhat
related. He showed that, if $\ct$ is the homotopy category
of spectra, and $H:\ct^{op}\ri\ab$ is a 
homological functor taking 
coproducts to products, then $H$ is representable. There are two
theorems here, one about homological
functors on $\ct^{op}$, and another about homological functors on the
subcategory $\{\ct^c\}^{op}$. And both theorems usually go under the
name Brown representability. 
Neither theorem is a special case of the other.
In the literature, one sometimes distinguishes them by
calling the theorem
about functors on $\{\ct^c\}^{op}$ ``Brown representability
for homology'', while the theorem about functors on $\ct^{op}$ goes
by the name ``Brown representability for cohomology''. 

The reason for this strange terminology is the following.
In many interesting cases, the category $\ct^c$ is self dual.
Thus functors ${\{\ct^c\}}^{op}\ri\ab$
are the same as functors $\ct^c\ri\ab$,
and these correspond to functors $\ct\ri\ab$
respecting coproducts. Thus ``Brown representability
for homology'' may be viewed as a theorem about
covariant homological functors $\ct\ri\ab$,
respecting coproducts.
\ermk

In hindsight, it seems natural to ask how these statements 
generalise to other triangulated categories, in particular,
the derived category of a ring. 
Surprisingly, questions of this sort were not asked until the
1980's. 

Even then, the first questions to be asked were: To what extent
can results about rings be generalised to homotopy theory.
The first to suggest that this might be a fruitful pursuit
was probably Waldhausen. Waldhausen proposed that techniques
from homological algebra---Hochschild homology and cohomology,
trace maps, and cyclic versions of these---should all be done
in the context of $E^\infty$ ring spectra. The work that followed,
by Goodwillie, B\"okstedt, Hsiang, Madsen and many others since,
showed how good the idea was.

The idea that translating results from homotopy theory to
derived categories could be worthwhile came later. The first
paper we are aware of is Hopkins'~\cite{Hopkins85}; in it, 
one has a derived category version of the nilpotence theorem. 
But it was really only in B\"okstedt--Neeman's~\cite{Bokstedt-Neeman93} 
that the first attempt was made, to use homotopy theoretic techniques
to solve standard problems on derived categories. In the 1990's,
we have seen explosive growth in the subject. In 
\cite{Neeman92A} and \cite{Neeman96}, Neeman applied
techniques coming from homotopy theory
to the study of, respectively, the localisation theorem
in $K$--theory and to Grothendieck duality. The articles by
Rickard~\cite{Rickard97}, Benson, Carlson and 
Rickard~\cite{Benson-Carlson-Rickard95}, \cite{Benson-Carlson-Rickard96},
\cite{Benson-Carlson-Rickard97}, 
Benson and Krause~\cite{Benson-Krause99}, 
Krause~\cite{Krause99a},
and Benson and Gnacadja~\cite{Benson-Gnacadja},
give beautiful applications to group cohomology. 
Keller~\cite{Keller98}, \cite{Keller99}, \cite{Keller98a}
applies the techniques to the study of cyclic homology.
And Voevodsky \cite{Voevodsky95},  \cite{Voevodsky96} 
and \cite{Voevodsky98},
Suslin--Voevodsky~\cite{Suslin-Voevodsky96}
and Morel~\cite{Morel98} and \cite{Morel99},
have produced a string of results,
which apply homotopy theory to the study of motives.

Along with the applications, came the study of the degree 
to which the theorems extend. Homotopy theorists, over a
period of 30 years, developed certain tools to handle 
the category of spectra. It became interesting to know
which parts of these tools work, in the new and greater
generality. This has also led to a series of papers. 
Hovey, Palmieri and Strickland~\cite{Hovey-Palmieri-Strickland}
set up a convenient axiomatic formalism. Without 
going into detail, we remind the reader of the work of
Beligiannis~\cite{Beligiannis99}, Christensen~\cite{Christensen96},
Christensen--Strickland~\cite{Christensen-Strickland96}, 
Franke~\cite{Franke97}, Keller~\cite{Keller94a}, 
Krause~\cite{Krause98}, \cite{Krause99}, \cite{Krause99a},
Krause and Reichenbach~\cite{Krause-Reichenbach98},
and Neeman~\cite{Neeman4}, \cite{Neeman96A}
and \cite{Neeman99}.

This skimpy historical survey was intended to explain
why people have studied whether
Brown representability generalises to derived categories. 
As we mentioned in Remark~\ref{R00.2}, the term Brown
representability is used to cover two theorems. 
Brown representability for cohomology is a characterisation
of representable functors $\ct^{op}\ri\ab$, while
Brown representability for homology is a more complicated
statement about functors ${\{\ct^c\}}^{op}\ri\ab$. Of these, the
generalisation of Brown representability for cohomology
is very well understood. The best and most recent results
were obtained independently by Franke~\cite{Franke97}
and Neeman~\cite{Neeman99}, and one of the remarkable
aspects of their theorems, is that they prove new results
even in homotopy theory. The theorems tell us, that
Brown representability for cohomology generalises to
the categories of $E$--acyclic spectra and $E$--local
spectra, for any homology theory $E$.

This paper addresses the less well understood problem, of
Brown representability for homology. In the remainder of
the Introduction, we will do two things. First, we will
go through the history of this problem in detail, explaining
what was already known. Then, we will outline the counterexamples
and results obtained in this article. But before we start, we need
to establish some notation.

\ntn{N00.3}
All rings will be associative, with unit. All $R$--modules
will be right, unitary modules. The
ring $R$ is called {\em hereditary} if its
global dimension is at most $1$. The triangulated category
$\ct=D(R)$ will be the unbounded derived category of
right $R$--modules. The category $\ct^c$ is, as above, the
full subcategory of compact objects in $\ct$.

We will denote the category of right $R$--modules by the
symbol $\Mod{R}$. The subcategory of finitely presented
$R$--modules will be denoted $\mod{R}$. The category of all
additive functors 
${\{\ct^c\}}^{op}\ri\ab$ will be denoted $\mr$,
while the category of all additive functors 
${\{\mod{R}\}}^{op}\ri\ab$ will bear the name
$\Md(\mod{R})$.

When speaking of objects of the category $\mr$, that is,
of functors
${\{\ct^c\}}^{op}\ri\ab$, we frequently wish to 
single out the ones that are homological, that is,
take triangles to long exact sequences. We will 
feel free to interchangeably use the adjectives
``homological'', ``exact'' or ``flat''. We remind the
reader that an object of \mr is exact if and
only if it is a filtered colimit of representable functors. 
Furthermore, the representable functors are projective.
(We use
the term ``representable'' to mean functors of the form $\y{C}$,
with $C$ compact. In the literature, people
sometimes call all functors $\y X$ representable.)

We also need to recall the notion of purity for 
$R$--modules. A short exact sequence of $R$--modules
\[
\CD
0 @>>> A @>>> B @>>> C @>>> 0
\endCD
\]
is called {\em pure exact}, if it remains exact when
tensored with an arbitrary left $R$--module. Equivalently,
it is a pure exact sequence if, for every finitely presented module
$P$, the functor $\text{Hom}(P,-)$ takes it to an exact sequence
\[
\CD
0 @>>> \text{Hom}(P,A) @>>> \text{Hom}(P,B) @>>> \text{Hom}(P,C) @>>> 0.
\endCD
\]
An $R$--module
$P$ is called {\em pure projective}, if the functor $\Hom(P,-)$
takes pure exact sequences to exact sequences. 
A module $P$ is pure projective if and only if it is a summand of 
a coproduct of finitely presented modules.  The {\em pure
projective dimension} of an $R$--module $M$ is defined to
be the length of its shortest pure resolution by pure
projectives.

A module $I$ is called {\em pure injective,} if the functor
$\Hom(-,I)$ takes pure exact sequences to exact sequences.
The pure injective dimension of a module $I$ is the length
of the shortest pure resolution by pure injectives. The
pure global dimension of $R$, denoted $\pgldim{R}$, is
the supremum over all $M$, of the pure projective
dimension of $M$. This equals the supremum of the
pure injective dimensions. We refer the reader to
\cite{Jensen-Lenzing} for a more thorough discussion,
with proofs.

Finally, recall our shorthand: for $X\in\ct$, we write
$\y X$ for the exact=homological=flat functor
${\ct(-,X)\big|}_{\ct^c}$. It is also convenient to make
a definition which is not so standard:
\entn

\dfn{D00.4}{\bf (Beligiannis~\cite{Beligiannis99})}\ \ 
The pure global dimension of $\ct$, denoted
$\pgldim{\ct}$, is defined to be the supremum,
over all $X\in\ct$, of the projective dimension
in $\mr$ of the object $\y X$.
\edfn

The following proposition will be useful.

\pro{P00.5}{\rm (Beligiannis~\cite[Prop.~11.2]{Beligiannis99}.
The proof is based on an idea by Jensen, 
which appeared in a paper by
Simson~\cite[Thm.~2.7]{Simson77}.)}\newline
The pure global dimension of $\ct$ is also
the supremum over all homological=exact functors
$F$, of the projective dimension of $F$. Note that,
as we will discover in this article, there can be more $F$'s
than $\y X$'s.
\epro

Let $\ct$ be a triangulated category with coproducts,
and $\ct^c\subset \ct$ the full subcategory of compact
objects. We adopt the following notation: 
 
\begin{description}
\item [{[BRO]}] The category $\ct$ satisfies [BRO] if
every exact functor ${\{\ct^{c}\}}^{op} \ra \ab$ is of
the form $\y X$, for some $X\in\ct$.
\item [{[BRM]}] The category $\ct$ satisfies [BRM] if
every natural transformation $\y X \ra \y X'$ is induced
by a map $X \ra X'$.
\end{description}

\nin
The theorem  of Adams (see \ref{T00.1}) says, that if $\ct$ is
the homotopy category of spectra, then both [BRO] and
[BRM] hold in $\ct$. In \cite{Neeman4}, Neeman found a necessary
and sufficient condition for this to generalise, to arbitrary
compactly generated $\ct$'s. For this article, in the
statements that follow, assume $\ct=D(R)$ is the derived
category of a ring $R$.

\thm{th:neeman}\cite{Neeman4} The following are equivalent:
\begin{roenumerate}
\item
Both [BRM] and [BRO] hold in $\ct$
\item $\pgldim{\ct}\leq1$.
\end{roenumerate}
\ethm

\nin
The direction (i)$\Longrightarrow$(ii) was also observed 
in~\cite{Christensen-Strickland96}. Beligiannis, using 
his Proposition~\ref{P00.5} above, recently showed:

\thm{th:beli}\cite[Theorem 11.8]{Beligiannis99}
[BRM]$\Longrightarrow$[BRO].
\ethm
\nin
  Neeman~\cite{Neeman4} also showed that when $R$
is countable, [BRM] (and therefore also [BRO]) holds.

Keller produced the first example, where [BRM] fails. 
It may be found in Neeman's~\cite{Neeman4}. The example
hinges on the following observation. If [BRM] holds, then
by Theorem~\ref{th:neeman}, we have $\pgldim{\ct}\leq1$.
That is, for any object $X\in\ct$, $\y X$ has
projective dimension at most $1$. If $R$ is a noetherian ring,
this means that the cohomology modules $H^iY$ have
pure projective dimension at most $1$. For a counterexample,
one needs only produce an object $Y\in\ct=D(R)$, so that
its cohomology is of pure projective dimension greater than $1$.

The most recent progress preceding this article is a theorem of
Beligiannis:

\thm{T00.9}\cite[Remark 11.12]{Beligiannis99}
[BRO] holds, whenever $\pgldim{\ct}\leq2$.
\ethm

This leaves several obvious questions:
\bd
\item[Q1] What is the precise relation between the pure global
dimension of $R$, denoted $\pgldim R$, and the pure
global dimension of $\ct$, denoted $\pgldim{\ct}$?
\item[Q2] Just how closely are the two related to [BRM] and [BRO]?
\item[Q3] Does [BRO] hold in general?
\ed

In this article, we make progress on these questions.
Regarding Q1, we prove that for many rings $\pgldim R \leq \pgldim{\ct}$,
and that for hereditary rings this is an equality.
Then we give examples to show that in general the inequality can be
strict.

Regarding Q2, we give a precise relationship between pure global
dimension, [BRO] and [BRM] for hereditary rings.
Then we give examples to show that in general no such simple 
relationship holds.  At the same time we show that [BRO] can
fail, answering Q3.  For example, it fails for $R = k[x,y]$
when $k$ has cardinality at least $\aleph_{3}$ (Example~\ref{E2.11}).

Here is a more detailed overview of these results.
We begin with an easy proposition giving our positive
results about Q1.  It is followed by a description of
our counterexamples.
We end with our positive results about Q2.

\medskip

\nin
{\bf Proposition~\ref{P1.2}}\ \ {\em
\begin{roenumerate}
\item 
Suppose that $R$ is a coherent ring, and that all finitely
presented $R$--modules are of finite projective
dimension. (This hypothesis holds when $R$ is noetherian
of finite global dimension.)  Then we have
\[
\pgldim{R} \leq \pgldim{D(R)}.
\]
\item Suppose that $R$ is hereditary. Then we have
\[
\pgldim{R} = \pgldim{D(R)}.
\]
\end{roenumerate}
}

\medskip

Weaker versions of this proposition were known before,
and the inequality was after all at the basis of Keller's
counterexample to [BRM]. The really new result we show 
in this article is that,
for some $R$, the inequality can be strict; Example~\ref{E1.3}
gives such an $R$. The idea of the counterexample is to
produce two rings $R$ and $S$, of different pure global dimensions,
but with $D(R)\cong D(S)$.
Then $\pgldim D(R)=\pgldim D(S)$
must be at least the maximum, and strictly bigger than the
minimum, of $\pgldim R$ and $\pgldim S$.
These rings, due to Thomas Bruestle,
are finite-dimensional non-commutative $k$-algebras
described by means of quivers.

Even more, we show that in general the answer
to Q3 is negative: [BRO] can fail. 
It fails for the rings $R$ and $S$ mentioned above when the cardinality
of $k$ is at least $\aleph_{2}$,
for the ring $k[x,y]$ when $|k| \geq \aleph_{3}$,
and also for the ring $T = k\angles{X,Y}$ of polynomials in two
non-commuting variables when $|k| \geq \aleph_{2}$.
(In particular, since it is consistent with ZFC that $|\C| = \aleph_{3}$,
it is impossible to prove [BRO] using ZFC when $R = \C[x,y]$.)
The proof that these are counterexamples is
presented in Section~\ref{S2}. Our method is to find an
exact sequence
\[
\CD
0 @>>> \y A @>>> F @>>> \y B @>>> 0
\endCD
\]
in \mr,
and show that $F$ is not isomorphic to $\y Y$ for any $Y$. 
The idea is to study the extension group $\Ext^1(\y B,\y A)$.
We get a handle on this group using several spectral sequences.
The precise statement of our theorem is:

\medskip

\nin
{\bf Theorem~\ref{T2.10}.}\ \ 
{\em
Let $R$ be an associative ring. Assume that $R$ is coherent,
and that every finitely presented $R$--module
has a finite projective resolution.
Suppose there exists an $R$--module $N$ so that
\[
\text{\rm pure inj dim}(N) - \text{\rm inj dim}(N)\geq2.
\]
Then [BRO] fails for in $D(R)$. This means that
there exists a homological functor
\mbox{$F:{\{\ct^c\}}^{op}\ri\ab$}, which is not the restriction
of any representable. That is, there exists no $Y$
with $\y Y=F$. }

\medskip

What is mysterious here, is that given a homological $F$,
we cannot directly tell whether it is of the form $\y X$.
We have no criterion to distinguish $\y X$'s from other
homological functors. In fact, Beligiannis' Proposition~\ref{P00.5}
tells us, that given any homological $F$, there exists a
$\y X$ of projective dimension greater than or equal to that
of $F$; projective dimension will not distinguish $\y X$'s
from other homological functors. What we do amounts to finding a
trick,
to get around this problem.

For general rings, this is all we can say. We can give a refinement of
the results for hereditary rings $R$; recall that $R$ is {\em hereditary}
if its global dimension is $\leq 1$. Examples of
hereditary rings are commutative principal ideal domains,
and non-commutative polynomial rings.

\medskip

\nin
{\bf Theorem~\ref{th:hereditary}.}\ \ 
{\em
Let $R$ be a hereditary ring.
Then
\begin{roenumerate}
\item {[BRM]} holds in $\ct$ if and only if the pure global dimension of $R$ 
is at most $1$; and
\item {[BRO]} holds in $\ct$ if and only if the pure global dimension of $R$ 
is at most $2$.
\end{roenumerate}
}

\medskip

We conclude the paper with the observation (Lemma~\ref{L2.12}) that
any counterexample to [BRO] must take values in infinite-dimensional
vector spaces.

\smallskip

\emph{Acknowledgements.}  The authors would like to thank 
Apostolos Beligiannis, Thomas Br\"ustle,
Henning Krause and
Michel Van den Bergh
for helpful conversations.
The first and second authors thank the third author, and the Centre 
for Mathematics and its Applications at the Australian National
University, for providing a friendly and productive setting while
this work was carried out.

\section{Pure global dimension: module categories versus derived
  categories}
\label{S1}

Let $R$ be an associative ring. 
We denote by $\ct$ the unbounded derived category
$D(R)$ of the category of (right) $R$--modules, and by $\ct^c$ the full
subcategory of compact objects. Recall that a complex is a compact
object of $\ct$ iff it is quasi-isomorphic to a bounded complex of
finitely generated projective $R$--modules. Here and elsewhere, we
identify the category $\Mod{R}$ of $R$--modules with the subcategory of
$\ct$ consisting of complexes concentrated in degree $0$.

\lem{L1.0} The following are equivalent:
\begin{roenumerate}
\item $R$ is coherent and each finitely presented $R$--module
is of finite projective dimension.
\item Each finitely presented $R$--module is compact when
viewed as an object of $D(R)$.
\item A complex $X$ is compact iff each $H^{n} X$ is finitely presented and
$H^nX\cong0$ for all but finitely many $n$.
\end{roenumerate}
\elem

\rmk{1.0} In particular, the conditions of the lemma are
satisfied if $R$ is noetherian and of finite global 
dimension.  They are also satisfied by any hereditary ring,
that is, any ring of global dimension at most 1.
\ermk

\prf 
We will prove
(i)$\Longleftrightarrow$(ii), and then that
(i)+(ii)$\Longleftrightarrow$(iii).
But first, we remind the reader 
that a ring is coherent iff the kernel of every map
between finitely generated projective modules is finitely
presented. We will also use the easy fact that
a module is a compact object of $D(R)$ iff it
admits a finite resolution by finitely generated projective
objects.

 Assume (i)  holds. Let $M$ be a finitely
 presented module. Since $R$ is coherent, $M$ admits a
 resolution by finitely generated projective modules.
 Since $M$ is of finite projective dimension, this
 resolution may be chosen to be finite. So 
 $M$ is compact in $D(R)$. That is, (ii) follows.

Suppose that (ii) holds. 
Then each finitely presented
module admits a finite resolution by finitely generated
projectives, and so in particular has finite projective dimension.
Now let $K$ be the kernel of a map $f:P_1 \ri P_0$ between finitely
generated projectives. Let $C$ be the cokernel
of $f$. In $D(R)$, we have the canonical triangle
\[
\Sigma K \ra P \ra C \ra \Sigma^2 K,
\]
where $P$ is the complex $P_1 \ra P_0$. By assumption,
$P$ and $C$ are compact. Hence $K$ is compact.
So it admits a finite resolution by finitely 
generated projective objects. In particular,
it is finitely presented. Thus $R$ is coherent; (i) holds.

 Thus far, we have proved (i)$\Longleftrightarrow$(ii).
 Assume these equivalent conditions hold; we wish to prove
 (iii). Let $X$ be a compact object in $D(R)$. It is 
 isomorphic to a
 finite complex of finitely generated projective modules. By
 (i), $R$ is coherent; hence $H^nX$ is finitely presented for 
 all $n$. And since the complex $X$ is finite,
 $H^nX\cong0$ for all but finitely many $n$.

 Suppose now that $H^nX$ is finitely presented for 
 all $n$, and that
 $H^nX\cong0$ for all but finitely many $n$.
 The $t$--structure on $D(R)$ gives us triangles
 \[
 \CD
 X^{\leq n} @>>> X @>>> X^{>n} @>>> \Sigma X^{\leq n} 
 \endCD
 \]
 and these allow us to assemble $X$ from its homology.
 Now $H^nX$ is finitely presented for 
 all $n$, and by (ii) it is compact. This forces $X$, an
 iterated extension of compact objects, to also be compact.
 We conclude that (iii) holds.

Finally, (iii)$\Longrightarrow$(ii) is immediate.
\eprf

 Recall that the functor
 $\y{ } : \ct \rightarrow \Modtc$
 sends an object $X\in\ct$ to the functor
\[
 \y{X}=\,\ct\,(-,X)|_{\ct^c}.
 \]
 For $i\in\zz$ and $F\in\Modtc$, we define the
 {\em $i^{th}$ homology of $F$} by
\[
H^i F = F(\Sigma^{-i} R).
\]
The functor $H^i: \Modtc \ri \Mod{R}$ extends the homology functor on $\ct$
in the sense that we have a canonical isomorphism $H^i \circ \y{ } =
H^i$.
 
An object $G$ in the category $\Modtc$ is called {\em finitely
presented,} if there exists an exact sequence
\[
 \y X \ra \y Y \ra G \ra 0.
\]
The full 
subcategory of all finitely presented objects in
$\Modtc$ is known 
to be an abelian category; see for example Freyd's~\cite{Freyd66A}.
As in the case of a module category, a sequence 
\[
0 \ra F_1 \ra F_2 \ra F_3 \ra 0
\]
of $\Modtc$ is called {\em pure exact} if the sequence
\[
0 \ra \Hom (G,F_1) \ra \Hom(G,F_2) \ra \Hom(G,F_3) \ra 0
\]
is exact for each finitely presented functor $G$. (In particular,
the sequence is then exact.) 

\lem{L1.1}
Suppose that the conditions of Lemma~\ref{L1.0} hold.
\begin{roenumerate}
\item The functor $\y{ }:\Mod{R}\ri \Modtc$ 
commutes with filtered colimits.
It takes pure projective $R$--modules
to projective objects of $\Modtc$. It transforms pure
exact sequences of $R$--modules into pure exact sequences in
$\Modtc$.                                                \label{it:y}
\item For each $i\in\zz$, the functor $H^i$ commutes with
filtered colimits. It takes projective
objects of $\Modtc$ to pure projective $R$--modules.
It transforms pure exact sequences of 
$\Modtc$ into pure exact sequences of $R$--modules.       \label{it:H}
\end{roenumerate}
\elem

\prf \iref{it:y} Let $M_\lambda$ be a filtered system of
$R$--modules. Clearly, if $P=\Sigma^i R$ for some $i\in\zz$, 
the canonical map
\[
\colim\, \ct\,(P, M_\lambda) \ri \ct\,(P, \colim\,M_\lambda)
\]
is bijective. Since both sides are cohomological functors
of $P$, this map is still bijective if $P$ is any compact
object of $\ct$, since $\ct^{c}$ is the thick subcategory generated by $R$.
This means that $\y{ }$ takes $\colim\,M_\lambda$
to $\colim\,\y{M_\lambda}$.

Each pure projective $R$--module is a direct factor of a
coproduct of finitely presented modules. Since the functor
$\y{ }$ commutes with coproducts, it is enough to show
that $\y{M}$ is projective if $M$ is finitely presented.
But in this case, $M$ is compact in $\ct$, by our assumption
on the ring $R$. So $\y{M}$ is projective since it is 
even representable.

Now let 
\[
0 \ri L \ri M \ri N \ri 0
\]
be a pure exact sequence of $R$--modules. Clearly, if $N$ is
finitely presented, the sequence splits. An arbitrary module
$N$ is a filtered colimit of finitely presented modules.
Thus the sequence is a filtered colimit of split sequences.
Since the functor $\y{ }$ commutes with filtered colimits,
the image of the sequence is also a filtered colimit of
split sequences. Thus it is pure.
\bigskip

\iref{it:H} By definition, the functor $H^i$ is evaluation
at $\Sigma^{-i} R$. Thus it commutes with colimits.
The projective objects of $\Modtc$ are direct factors
of coproducts of representable functors, and the functor
$H^i$ commutes with coproducts. So it is enough to show
that $H^i \y{P}=H^i P$ is pure projective for $P\in \ct^c$.
This is clear since $H^i P$ is finitely presented, by
our assumption on the ring $R$.

Let 
\[
0 \ra F_1 \ra F_2 \ra F_3 \ra 0
\]
be a pure exact sequence of $\Modtc$. Clearly if $F_3$ is
finitely presented, the sequence splits. In the general case,
$F_3$ is a filtered colimit of a system of finitely
presented functors. So the sequence is a filtered colimit
of split sequences. Since the functor
$H^i$ commutes with filtered colimits, this implies the
last assertion.
\eprf

The {\em pure global dimension} of the derived category
$D(R)=\ct$ is by definition~\cite{Beligiannis99} 
the supremum of the projective dimensions
of the functors $\y{X}$, $X\in\ct$. We write
$\pgldim{ } $ for `pure global dimension'.
Part (ii) of the following lemma is 
due to Beligiannis~\cite[Prop. 12.8]{Beligiannis99}.

\pro{P1.2} 
Suppose that the conditions of Lemma~\ref{L1.0} hold.
\begin{roenumerate}
\item Let $M$ be an $R$--module. Then the projective
dimension of $\y{M}$ equals the pure projective dimension of $M$.
Hence we have
\[
\pgldim{R} \leq \pgldim{D(R)}.
\]
\item Suppose that $R$ is hereditary. Then we have
\[
\pgldim{R} = \pgldim{D(R)}.
\]
\end{roenumerate}
\epro

\prf
(i) The first part of the preceding lemma shows that the functor
$\y{ }$ takes pure projective resolutions of a module $M$ to projective
resolutions of $\y{M}$. Hence the projective dimension of 
$\y{M}$ is no more than the pure projective dimension of $M$.
Conversely, let 
\[
\cdots \ri Q_1 \ri Q_0 \ri \y{M}\ri 0
\]
be a projective resolution of $\y{M}$. If $M$ is finitely presented,
then $\y M$ is projective, so
the resolution is nullhomotopic. An arbitrary $M$ is still a filtered
colimit of finitely presented modules. So for arbitrary $M$ the
resolution is a filtered colimit of nullhomotopic complexes.
Thus it is a pure exact sequence. By the second part
of the above lemma, its image under $H^0$ is a pure
projective resolution of $H^0 \y{M}= M$. Thus the pure projective
dimension of $M$ is no more than the projective dimension of $\y{M}$.

(ii) By part (i), it suffices to prove that $\pgldim{R} \geq \pgldim{D(R)}$.
Let $X\in D(R)$. Since $R$ is hereditary, the object $X$ is
isomorphic in $D(R)$ to the coproduct of the $\Sigma^{-i} H^i X$,
$i\in\zz$; see Lemma~6.7, on page 153 of
\cite{Neeman92}. Hence the projective dimension of 
$\y{X}$ is no greater than the supremum of the projective dimensions 
of the $\y{H^i X}$.
These are bounded by $\pgldim{R}$ thanks to part (i).
\eprf

\exm{E1.3} Let $k$ be a field and let $t$ be the cardinal such
that $\aleph_t = \max(|k|,\aleph_0)$.
So $t$ is $0$ if $k$ is finite or countable, $1$ if $k$ has the
smallest uncountable cardinality, etc.
Building on an example due to Th.~Bruestle
we will exhibit a $k$-algebra $R$ such that the inequality
\[
\pgldim R \leq \pgldim D(R)
\]
is strict. Our example
is based on the observation that there are algebras with 
equivalent derived categories but widely differing
pure global dimensions. More precisely, 
we will exhibit a finite-dimensional 
$k$-algebra $R$ with $\pgldim R=0$ such that
$D(R)$ is triangle equivalent to $D(S)$
for a finite-dimensional hereditary 
$k$-algebra $S$ whose pure global dimension is $t+1$
($\infty$ if $t$ is infinite).
Thus we have
\[
\pgldim{R} < \pgldim{S} = \pgldim{D(S)} = \pgldim{D(R)},
\]
where we have used part (ii) of the above proposition
for the first equality.

Thus Theorem~\ref{th:hereditary} implies that 
[BRM] fails for $D(R)$ when $t \geq 1$ and that
[BRO] fails for $D(R)$ when $t \geq 2$, 
even though $R$ has pure global dimension $0$.

The algebras $R$ and $S$ are due to Th.~Bruestle.
We will define them using the language of quivers
with relations (cf.~\cite{Ringel84}, \cite{GabrielRoiter92},
\cite{AuslanderReitenSmaloe95}). Here is all we need: 
A {\em quiver} is an oriented graph. It is
thus given by a set $Q_0$ of points, a set $Q_1$ of arrows,
and two maps $s,t: Q_1 \ri Q_0$ associating with each
arrow its source and its target. A simple example is the
quiver
\[
\vec{A}_{10}\; : \quad 1 \xrightarrow{\alpha_1} 2 \xrightarrow{\alpha_2} 3 \ri
\cdots \ri 8 \xrightarrow{\alpha_8} 9 \xrightarrow{\alpha_9} 10.
\]
A {\em path} in a
quiver $Q$ is a sequence
$(y|\beta_r| \beta_{r-1} | \cdots | \beta_1| x)$ of
composable arrows $\beta_i$ with $s(\beta_1)=x$,
$s(\beta_i)=t(\beta_{i-1})$, $2\leq i\leq r$, $t(\beta_r)=y$.
In particular, for each point $x\in Q_0$, we have
the {\em lazy path} $(x|x)$. It is neutral for the
obvious {\em composition} of paths. The {\em quiver algebra}
$kQ$ has as its basis all paths of $Q$. The product
of two basis elements equals the composition of the
two paths if they are composable and $0$ otherwise.
For example, the quiver algebra of $Q=\vec{A}_{10}$ is isomorphic
to the algebra of lower triangular $10\times 10$ matrices.

The construction of the quiver algebra $kQ$ is motivated
by the (easy) fact that the category of left $kQ$-modules
is equivalent to the category of all diagrams of
vector spaces of the shape given by $Q$. It is
not hard to show that each quiver algebra is
hereditary. It is finite-dimensional over $k$ iff the
quiver has no oriented cycles.

Gabriel \cite{Gabriel72} showed that the quiver algebra
of a finite quiver has only a finite number of 
$k$--finite-dimensional indecomposable modules 
(up to isomorphism) iff the underlying graph
of the quiver is a disjoint union of Dynkin diagrams
of type $A$, $D$, $E$.

The above example has underlying graph of Dynkin type $A_{10}$ and
thus its quiver algebra has only a finite number of
finite-dimensional indecomposable modules.

An ideal $I$ of a finite quiver $Q$ is {\em admissible}
if for some $N$ we have 
\[
(kQ_1)^N \subseteq I \subseteq (kQ_1)^2,
\]
where $(kQ_1)$ is the two-sided ideal generated by all paths of length
$1$.  A \emph{quiver $Q$ with relations $R$} is a quiver $Q$ with a
set $R$ of generators for an admissible ideal $I$ of $kQ$. The algebra
$kQ/I$ is then the \emph{algebra associated with $(Q,R)$}.  Its
category of left modules is equivalent to the category of diagrams of
vector spaces of shape $Q$ obeying the relations in $R$. The algebra
$kQ/I$ is finite-dimensional (since $I$ contains all paths of length
at least $N$), hence artinian and noetherian.  
By induction on the number of points one can show that
if the quiver $Q$ contains no oriented cycle, then the
algebra $kQ/I$ is of finite global dimension.

One can show that every finite-dimensional algebra
over an algebraically closed field is Morita equivalent
to the algebra associated with a quiver with relations
and that the quiver is unique (up to isomorphism).

Now we let $R$ be the finite-dimensional
$k$-algebra associated with the above quiver $\vec{A}_{10}$
and the relation $\alpha_8 \alpha_7 \cdots \alpha_1$
(no $\alpha_9$!). The algebra $R$ is a quotient of $k \vec{A}_{10}$ 
and thus it admits only a finite number of indecomposable finite-dimensional
modules. By a result of Auslander \cite{Auslander74} and Tachikawa
\cite{Tachikawa73}, this is equivalent to $\pgldim{R}=0$.

Let $S$ be the quiver algebra of the quiver 
\[
E : 
\begin{array}{lllllllll}
2 \ri & 3 \ri & 4 \ri & 5 \ri & 6 \ri & 7 \ri & 8 \ri & 9 \ri & 10 \\
      &       &       &       &       &       & \da   &       &    \\
      &       &       &       &       &       & 1.    &       &
\end{array}
\]
Thus $S$ is finite-dimensional over $k$ and hereditary. 
By Theorem 4.1 of Baer-Lenzing's \cite{BaerLenzing82}, 
we have $\pgldim S=t+1$ ($\infty$ if $t$ is infinite).

Finally, we need to show that $R$ and $S$ have equivalent
derived categories. Indeed, the algebra $R$ 
admits a tilting complex with endomorphism ring $S$
so that the equivalence follows from Rickard's
Morita theorem for derived categories \cite{Rickard89b}.
To describe the tilting complex, let $P_i = e_i R$
be the projective $R$--module associated with the
idempotent $e_i=(i|i)$ (the lazy path). 
It is easy to compute the morphism spaces between
these modules: Indeed, we have $\Hom(e_i R, e_j R) = e_j R e_i$
and this space identifies with the vector space on the
set of paths from $i$ to $j$ divided by the subspace
of linear combinations of paths lying in the ideal
of relations. For example, for $i\leq j$, the path from $i$ to $j$
yields a canonical morphism $P_i  \ra P_j$, which vanishes
iff $(i,j)=(1,9)$ or $(i,j)=(1,10)$. 
The tilting complex $T$
is now the sum of the complexes
\begin{align*}
& T_2=(P_1 \ri P_2),\quad T_3=(P_1 \ri P_3), \quad
\ldots\;\;, \quad T_8=(P_1 \ri P_8), \quad  \\
& T_1=(P_1 \ri 0),\quad T_9=(0 \ra P_9),\quad
T_{10}=(0 \ra P_{10}),
\end{align*}
where the first term of each complex is in degree $0$.
Using the description of the morphism spaces between
the $P_i$ it is not hard to check that, in the homotopy
category of right $R$-modules, we do have
$\Hom(T_i, T_j[l])=0$ for all $i,j$ and all $l\neq 0$,
and that the endomorphism ring of $T$ is indeed isomorphic
to $S$. For example, the canonical idempotent $(i|i)$
of the quiver $E$ corresponds to the idempotent
of $\End(T)$ arising from the identity of $T_i$ and
the arrow $8 \ri 9$ of $E$ corresponds to the
obvious morphism of complexes
\[
\begin{array}{rcl}
P_1 & \ri & P_8 \\
\da &     & \da \\
0   & \ra & P_9
\end{array}
\]
which is well-defined thanks to the relation
$\alpha_8 \alpha_7 \cdots \alpha_1$ that we imposed.
\eexm

\section{Failure of Brown representability}
\label{S2}

In this section, $R$ will be a ring satisfying the
equivalent conditions of Lemma~\ref{L1.0}.
In particular, all the theorems hold if
$R$ is a noetherian ring
of finite global dimension, or if $R$ is hereditary.
We begin by reminding ourselves of a standard spectral sequence.

\lem{L2.2}
Let $\ca$ be an abelian category satisfying AB5, and with
enough projectives. Suppose that $X$ and $Y$ are objects of 
$\ca$ and that 
$X=\colim \, X_\lambda$
expresses $X$ as a filtered colimit of objects $X_\lambda\in\ca$.
Then there is a spectral sequence, converging to
$\Ext^{i+j}(X,Y)$, whose $E^2$ term is
\[
\climi \Ext^j(X_\lambda,Y).
\]
\elem

\prf
There is a standard chain complex which computes the
derived functors of $\colim$.
Since the abelian category $\ca$ satisfies AB5, the derived
functors of filtered colimits vanish, and we deduce an
exact sequence in $\ca$
\[
\CD
\cdots @>>> \ds\bigoplus_{\lambda \sr \mu}X_\lambda
@>>> \ds\bigoplus_{\lambda}X_\lambda
@>>> X @>>> 0.
\endCD
\]
This gives us a resolution of $X$ in $\ca$, and the
spectral sequence is just the spectral sequence
of the functor $\Ext^{*}(-,Y)$ applied to this resolution.
\eprf

In the following, we write $\mod{R}$ for the category of finitely 
presented $R$--modules and $\Md(\mod{R})$ for the category of
contravariant additive functors from $\mod{R}$ to $\ab$.
The object 
\[
{\Mod{R}\big(-\,,\,M\big)\Big|}_{\mod{R}}
\]
of $\Md(\mod{R})$ will be denoted $\x M$.

\lem{L2.3}
Let $R$ be a ring, and let $M_{\lambda}$ be a filtered diagram of $R$--modules
with colimit $M$.
Then
\begin{roenumerate}
\item $\y{M} = \colim\ \y{M_{\lambda}}$ in \mr.
\item $\x{M} = \colim\ \x{M_{\lambda}}$ in $\Md(\mod{R})$.
\end{roenumerate}
\elem

\prf (i) was proved in Lemma \ref{L1.1} \iref{it:y}.
The second statement is more familiar in the equivalent form, which states that
$\Mod{R}(K,M) = \colim \Mod{R}(K,M_{\lambda})$ for any finitely presented $K$.
This is not hard to prove.
%
\eprf

\rmd{R2.4}
Let $R$ be a ring and let $M$ be an $R$--module.
Consider the filtered diagram of finitely presented modules $M_{\lambda }$
equipped with a map to $M$.
Then $M$ is the colimit of this diagram;
we already used this in the proof of Proposition~\ref{P1.2}(i).
This is the setting in which we will apply Lemma~\ref{L2.3}.
\ermd

The following lemma is well known; the proof may be found,
for example, in Theorem~2.8 of Simson's~\cite{Simson77}.
We include a sketch of the proof for the reader's convenience.

\lem{L2.4}
Let $R$ be a ring satisfying the conditions of Lemma~\ref{L1.0},
and let $M$ be an $R$--module. As mentioned in Remark~\ref{R2.4},
$M$ is the filtered colimit of all finitely presented
modules $M_\lambda$ mapping to $M$.
\begin{roenumerate}
\item \label{L2.4.1}
Let $F$ be an object of \mr. That is, $F$ is a functor
${\{\ct^c\}}^{op}\ri \ab$. Then the group 
$\Ext^i(\y{M},F)$
of extensions in \mr is isomorphic to $\climi F(M_\lambda)$.
\item \label{L2.4.2}
Let $F$ be an object of $\Md(\mod{R})$. That is, $F$ is a functor
$\{\mod{R}\}^{op}\ra \ab$. Then the group
$\Ext^i(\x{M},F)$
of extensions in $\Md(\mod{R})$ is isomorphic to $\climi F(M_\lambda)$. 
\end{roenumerate}
\elem

\prf
(i): By Lemma~\ref{L2.3}, $\y M$ is the colimit of $\y M_\lambda$ in \mr.  
Lemma~\ref{L2.2} then tells us that we get a spectral sequence with
$E^2$ term
\[
\climi\Ext^j(\y M_\lambda, F)
\]
converging to the group $\Ext^{i+j}(\y M, F)$ of extensions in \mr.
The functor $\y M_\lambda$ is representable, since by our hypothesis on 
$R$ the module $M_{\lambda}$ is compact.
Thus $\y M_{\lambda}$ is projective,
the $\Ext^j$ terms vanish unless $j=0$,
the spectral sequence collapses, and the desired isomorphism follows.

The proof of (ii) is similar. 
\eprf

\rmk{R2.5}
In part~\iref{L2.4.1} of Lemma~\ref{L2.4}, we computed the extensions of 
$\y M$ by $F$. This interests us most in the case where 
$F=\y{\Sigma^jN}$, with $N$ an $R$--module. 
In this case, the computation tells us that we have isomorphisms
\[
\Ext^i(\y{M},\y{\Sigma^jN}) =
\climi\ct(M_\lambda,\Sigma^jN) = 
\climi\Ext_R^j(M_\lambda,N).
\]

In part~\iref{L2.4.2} of Lemma~\ref{L2.4}, we computed the extensions of 
$\x M$ by $F$. This interests us most in the case where 
$F=\x N$, with $N$ an $R$--module. 
In this case, the computation tells us that we have an isomorphism
\[
\Ext^i(\x M,\x N) =
\climi\Hom_R(M_\lambda,N).
\]
Moreover the group $\Ext^{i}(\x M, \x N)$ 
of extensions in $\Md(\mod{R})$ can be identified with the group 
$\PExt^{i}(M,N)$; see \cite{Jensen-Lenzing}. We deduce that
\[
\PExt^{i}(M,N) =
\climi\Hom_R(M_\lambda,N).
\]
\ermk

\cor{C2.6}
If $M$ and $N$ are $R$--modules and $j>0$, then every map
$\y{\Sigma^j M} \ra \y{N}$
vanishes. 
Moreover, maps
$\y{M} \ra \y{N}$
are in one-to-one correspondence with maps of $R$--modules $M\ri N$.
\ecor

\prf
For $j>0$, we must show that any map
$\y{M} \ra \y{\Sigma^{-j}N}$
vanishes. But by Remark~\ref{R2.5}, the group of
such maps is
\[
\clim^{\!0}\,\Ext^{-j}_R(M_\lambda,N),
\]
which vanishes because there are no extensions
of negative degree. 

The group of maps 
$\y{M} \ra \y{N}$
is exactly
\[
\clim^{\!0}\,\Ext^{0}_R(M_\lambda,N),
\]
which is $\Hom_R(M,N)$.
\eprf

\lem{L2.7}
Let $F$ be an object in \mr, that is, a contravariant
additive
functor from $\ct^c$ to $\ab$. Suppose there exists
an integer $j>0$, $R$--modules $M$ and $N$, and a short
exact sequence in \mr
\[
\CD
0 @>>> \y{\Sigma^j N} @>\alpha>> F @>\beta>> \y M @>>> 0.
\endCD
\]
Then this sequence is unique up to isomorphism.
\elem

\prf
The integer $j$ and the modules $M$ and $N$ are clearly determined by
the homology of $F$.
In Corollary~\ref{C2.6} we saw that any map
$\y{\Sigma^j N} \ra \y M$
vanishes. Therefore, given any map
$\gamma : \y{\Sigma^j N} \ra F$,
the composite
 \[
\CD
\y{\Sigma^j N} @>\gamma>> F@>\beta>> \y M
\endCD
\]
vanishes, and hence $\gamma$ must factor through
$\alpha$. Dually, any map 
$F \ra \y M$
must factor through $\beta$. This shows that the given exact
sequence is unique.
\eprf

\lem{L2.8}
Let $F$ be an object of \mr, and suppose
there exists
an integer $j>0$, $R$--modules $M$ and $N$, and a short
exact sequence in \mr
\[
\CD
0 @>>> \y{\Sigma^j N} @>\alpha>> F @>\beta>> \y M @>>> 0.
\endCD
\]
The functor $F$ will be of the form $\y{Y}$ if and only
if the short exact sequence comes from a triangle. That is,
if and only if there exists a triangle in $\ct$
\[
\CD
\Sigma^j N @>>> Y @>>> M @>\partial>> \Sigma^{j+1}N
\endCD
\]
with $\partial$ a phantom map, so that the sequence
\[
\CD
0 @>>> \y{\Sigma^j N} @>\alpha>> F @>\beta>> \y M @>>> 0
\endCD
\]
is obtained by restricting the representable
functors to $\ct^c$.
\elem

We remind the reader that a map $W \ra X$ in $\ct$ is called
\df{phantom} if the composite $C \ra W \ra X$ is zero for each
compact object $C$ and each map $C \ra W$.

\prf
The implication $\Longleftarrow$ is trivial. If the
triangle exists and is isomorphic to
the short exact sequence of functors on
$\ct^c$, then $F$ is the restriction
of a representable functor on $\ct$. We wish to
prove $\Longrightarrow$. We suppose therefore
that the short exact sequence of functors is given, and that 
$F$ is the restriction of a representable. We want 
to produce a triangle.

The short exact sequence
\[
\CD
0 @>>> \y{\Sigma^j N} @>\alpha>> F @>\beta>> \y M @>>> 0
\endCD
\]
permits us easily to compute $F(\Sigma^nR)$, for all $n\in\zz$.
We have
\[
F(\Sigma^nR)\eq\left\{
\begin{array}{lcl}
M &\qquad& \mbox{if }n=0 \\
N &      & \mbox{if }n=j \\
0 &      & \mbox{otherwise. }
\end{array}
\right.
\] 
But if $F=\y Y$, then $F(\Sigma^nR)=H^{-n}(Y)$. The
above computes for us the cohomology of $Y$, as an object
in $D(R)=\ct$.

There is a $t$-structure truncation on $D(R)$, giving
a triangle
\[
\CD
Y^{\leq-1} @>>> Y @>>> Y^{\geq0} @>\partial>> \Sigma Y^{\leq-1} ,
\endCD
\]
and our homology computation shows that $Y^{\leq-1}$ and
$Y^{\geq0}$ each have only one non-zero cohomology group.
The triangle is therefore of the form
\[
\CD
\Sigma^j N @>>> Y @>>> M @>\partial>> \Sigma^{j+1}N.
\endCD
\]
We deduce an exact sequence
\[
\CD
\y{\Sigma^j N} @>>> \y Y @>>> \y M .
\endCD
\]
Now recall that $\y Y=F$, and that by the proof of
Lemma~\ref{L2.7},
any map $\y{\Sigma^j N}\ri F$ factors through
$\alpha$, and any map $F\ri\y M$ factors
through $\beta$. The exact sequence coming from the
triangle therefore factors through
\[
\CD
  @.   \y{\Sigma^j N} @.              @.             @.  \\
@.     @VfVV                     @.             @.        @. \\
0 @>>> \y{\Sigma^j N} @>\alpha>> F @>\beta>> \y M @>>> 0.\\
@.     @.                     @.             @VVgV        @. \\
  @.    @.              @.            \y M  @.  
\endCD
\]
By Corollary~\ref{C2.6}, the morphisms $f$ and $g$ in the 
diagram above come from maps of modules $N\ri N$ and $M\ri M$.
Evaluating the functors at $R$ and $\Sigma^jR$, 
we compute that both $f$ and
$g$ are isomorphisms. Hence the triangle gives rise to the
short exact sequence of functors, and $\partial$ must be
a phantom map.
\eprf

Next comes a spectral sequence argument. To help the reader,
we will first do the easy, baby case.

\pro{P2.9}
Let $R$ be a ring satisfying the conditions of Lemma~\ref{L1.0}.
Let $N$ be an $R$--module with injective dimension at most $1$
and pure injective dimension at least $3$.
Then in \mr there exists a homological functor
$F:{\{\ct^c\}}^{op}\ri\ab$ which is not the restriction
of any representable. That is, there exists no $Y$ with $\y Y=F$. 
\epro

\exm{E2.9.3}
Let $k$ be a field and $R$ the algebra of the quiver $E$
of Example~\ref{E1.3} (we called it $S$ there). 
Then $R$ is finite-dimensional over $k$ 
and hereditary, since it is the quiver algebra of a finite quiver. 
So all $R$--modules are of injective dimension at most $1$. 
Assume that $k$ is infinite of cardinality $\aleph_t$. 
Then by \cite{BaerLenzing82}, the pure global dimension of $R$ equals $t+1$
($\infty$ if $t$ is infinite). 
Thus when $t \geq 2$ there does exist an $R$--module satisfying the
assumptions of the proposition. 

Similarly, the ring $k\angles{X,Y}$ of polynomials in two 
non-commuting variables is an example when $t \geq 2$.
\eexm

To obtain examples where $R$ is commutative, we
will need to use Theorem~\ref{T2.10}, which is 
a refined version of the above proposition.

\prf
Because $N$ is of pure injective dimension at least $3$, there
exists a module $M$ and integer $n\geq3$, so that
$\PExt^n(M,N)\neq0$. If $n>3$, choose a pure exact
sequence 
\[
\CD
0 @>>> M' @>>> P @>>> M @>>> 0,
\endCD
\] 
with $P$ pure projective. Then $\PExt^n(M,N)=\PExt^{n-1}(M',N)$.
By a sequence of such dimension shifts, we may find an $M$
so that 
\[
\PExt^3(M,N)\neq0.
\]

By Remark~\ref{R2.4}, we may express $M$ as a filtered colimit
of finitely presented modules $M_\lambda$. By Lemma~\ref{L2.2},
applied this time to the category of $R$--modules, there is
a spectral sequence with $E^2$ term
$\climi\Ext^j_R(M_\lambda,N)$
converging to $\Ext^{i+j}_R(M,N)$.
We will now compute in this spectral sequence.

In Remark~\ref{R2.5}, we computed that 
\[
\clim^{\!3}\,\Ext^0_R(M_\lambda,N)=\PExt^3(M,N),
\]
and by the above, this does not vanish. On the other hand,
we know that $\Ext^3_R(M,N)=0$, since by hypothesis $N$ is of
injective dimension at most $1$. It follows that one
of the differentials in the spectral sequence into
the term
\[
\clim^3\Ext^0_R(M_\lambda,N)
\]
must be non-zero.

But there are only two differentials into this term,
one from $\clim^{\!1}\,\Ext^1$ and one from $\clim^{\!0}\,\Ext^2$.
The latter vanishes,
since by hypothesis $N$ is of injective dimension at most $1$.
It follows that
\[
\clim^{\!1}\,\Ext^1_R(M_\lambda,N)\neq0.
\]
But in Lemma~\ref{L2.4} we showed that this is the group
of extensions, in \mr,
\[
\CD
0 @>>> \y{\Sigma N} @>>> F @>>> \y M @>>> 0.
\endCD
\]
The group does not vanish so we may choose a non-trivial extension. 
Since $F$ is the extension of two homological functors, $F$ must be
homological. Now we will show that $F$ cannot be isomorphic
to a functor $\y Y$.

Lemma~\ref{L2.8} tells us that if $F$ is isomorphic
to $\y Y$, then there is a triangle in $\ct$
\[
\CD
\Sigma N @>>> Y @>>> M @>\partial>> \Sigma^2N
\endCD
\]
so that the exact sequence of functors above is isomorphic
to the one obtained from the triangle. But the map
$\partial:M\ri \Sigma^2N$ is an element of 
\[
\Ext^2(M,N)=0,
\]
and therefore the triangle splits. The exact sequence
of functors is not split, and we conclude that $F$ cannot
be isomorphic to any $\y Y$.
\eprf

The next Theorem is the more macho
computation with the same
spectral sequence.

\thm{T2.10}
Let $R$ be a ring satisfying the conditions of Lemma~\ref{L1.0}.
Suppose there exists an $R$--module $N$ so that
\[
\text{\rm pure inj dim}(N) - \text{\rm inj dim}(N)\geq2.
\]
Then [BRO] fails for in $D(R)$. This means that
there exists a homological functor
\mbox{$F:{\{\ct^c\}}^{op}\ri\ab$} which is not the restriction
of any representable. That is, there exists no $Y$
with $\y Y=F$. 
\ethm

\prf
Let $N$ be a module satisfying the hypotheses.
Let $n=\text{\rm inj dim}(N)$. Then
$\text{\rm pure inj dim}(N)\geq n+2$.
As in the proof of Proposition~\ref{P2.9},
we may choose a module $M$ with
$\PExt^{n+2}(M,N) \neq 0.$
We may also express
$M$ as a filtered colimit of finitely presented modules
$M_\lambda$.

Lemma~\ref{L2.2} gives us a spectral sequence, whose $E^2$ term
is
\[
\climi\Ext^j_R(M_\lambda,N),
\]
which converges to $\Ext^{i+j}_R(M,N)$. Once again, we have
that
\[
\clim^{\!n+2}\,\Ext^0_R(M_\lambda,N)=\PExt^{n+2}(M,N),
\]
and this does not vanish, by the choice of $M$.
But $\Ext^{n+2}_R(M,N)=0$, since $N$
is of injective dimension at most $n$, so there must be
a non-zero differential into the term
\[
\clim^{\!n+2}\,\Ext^0_R(M_\lambda,N).
\]

Now observe that 
\[
\clim^{\!0}\,\Ext^{n+1}_R(M_\lambda,N)=0,
\]
since $N$ is of injective dimension at most $n$. It follows that
for some $i$ with $1\leq i\leq n$, there is a non-zero
differential in the spectral sequence, from
\[
\climi\Ext^{n+1-i}_R(M_\lambda,N)
\]
to the term $\clim^{\!n+2}\,\Ext^0_R(M_\lambda,N)\neq0$.

Now recall the construction of our spectral sequence, from
Lemma~\ref{L2.2}. Since $M$ is the filtered colimit of
$M_\lambda$, there is an exact resolution of $M$
\[
\CD
\cdots @>>> \ds\bigoplus_{\lambda \sr \mu}M_\lambda
@>>> \ds\bigoplus_{\lambda}M_\lambda
@>>> M @>>> 0.
\endCD
\]
This resolution is a pure exact resolution by pure projectives.
(It is pure exact because it remains exact in the category
\Md(\mod{R}). And direct sums of finitely presented
modules $M_\lambda$ are pure projective.) 
By Lemma~\ref{L1.1}, it becomes an exact
resolution by projectives in the category \mr.

To simplify the notation, we will write the above resolution as
\[
\CD
\cdots @>>> P_2
@>>> P_1 @>>> P_0
@>>> M @>>> 0.
\endCD
\]
Let $K_i$ stand for the image of the map $P_{i} \ra P_{i-1}$.
In Lemma~\ref{L2.4} we showed that 
\[
\climi\Ext^{n+1-i}_R(M_\lambda,N)
\]
is the group of extensions
\[
\Ext^i(\y M,\y {\Sigma^{n+1-i}N}).
\]
But since the pure exact sequence
\[
\CD
0@>>>K_{i-1} @>>> P_{i-2}
@>>> \cdots @>>> P_0
@>>> M @>>> 0
\endCD
\]
remains exact in \mr, and the middle modules map to
projectives in \mr,
we deduce that the above extension group
is isomorphic to
\[
\Ext^1(\y {K_{i-1}}, \y {\Sigma^{n+1-i}N}).
\]
In other words, an element of the group
\[
\climi\Ext^{n+1-i}_R(M_\lambda,N)
\]
may be thought of as a short exact sequence in
\mr
\[
\CD
0 @>>> \y {\Sigma^{n+1-i}N}
@>>> F @>>> \y{K_{i-1}} @>>> 0.
\endCD
\]
We know that in the spectral sequence, for some $1\leq i\leq n$,
there is a non-zero differential 
\[
\CD
\climi\Ext^{n+1-i}_R(M_\lambda,N) \quad \supset \quad E \quad
@>\gamma>>
\quad \clim^{\!n+2}\,\Ext^{0}_R(M_\lambda,N)
\endCD,
\]
for a subgroup $E \subset \climi\Ext^{n+1-i}_R(M_\lambda,N)$.
What we will now show is that, if $\gamma(x)\neq0$, then
$x$ corresponds to an exact sequence
\[
\CD
0 @>>> \y {\Sigma^{n+1-i}N}
@>>> F @>>> \y{K_{i-1}} @>>> 0
\endCD
\]
where $F$ is not isomorphic to any $\y Y$.
Expressing the same thing slightly differently,
we will show that if $x\in\climi\Ext^{n+1-i}_R(M_\lambda,N)$
comes from an exact sequence of functors with $F=\y Y$,
then $\gamma(x)=0$.

Suppose therefore that we are given a short exact
sequence in \mr
\[
\CD
0 @>>> \y {\Sigma^{n+1-i}N}
@>>> \y Y @>>> \y{K_{i-1}} @>>> 0.
\endCD
\]
We need to show that in the spectral sequence, the differential $\gamma$
annihilates $x$. By Lemma~\ref{L2.8}, the exact sequence
of functors comes from a triangle
\[
\CD
\Sigma^{n+1-i}N @>>> Y @>>> K_{i-1} @>\partial>> \Sigma^{n+2-i}N
\endCD
\]
with $\partial$ a phantom map. 
 From the definition of the modules $K_i$, 
we have a pure exact sequence of $R$--modules
\[
\CD
0 @>>> K_i @>>> P_{i-1} @>>> K_{i-1} @>>> 0.
\endCD
\]
This exact sequence gives a triangle in $\ct=D(R)$. 
The fact that $\partial:K_{i-1}\ri \Sigma^{n+2-i}N$
is phantom tells us that the composite
\[
\CD
P_{i-1} @>>> K_{i-1} @>\partial>> \Sigma^{n+2-i}N
\endCD
\]
must vanish, since $P_{i-1}$ is a coproduct of compact
objects. But then the map $\partial$ must factor
as 
\[
\CD
K_{i-1} @>>> \Sigma K_i @>>> \Sigma^{n+2-i}N.
\endCD
\]
Thus if an element $x\in\climi\Ext^{n+1-i}_R(M_\lambda,N)$
comes from a short exact sequence
\[
\CD
0 @>>> \y {\Sigma^{n+1-i}N}
@>>> F @>>> \y{K_{i-1}} @>>> 0
\endCD
\]
with $F\simeq\y Y$, then $Y$ is determined by a class
\[
y\in\Ext^{n+1-i}_R(K_i,N).
\]

In conclusion, we deduce the following. Let us define 
$K_0=0$. We have a map of chain complexes
\[
\CD
\cdots @>>>P_i @>>> P_{i-1}
@>>> \cdots @>>> P_1
@>>> P_0 @>>> 0 \\
@. @VVV @VVV @. @VVV @VVV @. \\
\cdots @>0>>K_i @>0>> K_{i-1}
@>0>> \cdots @>0>> K_1
@>0>> K_0 @>>> 0 \period
\endCD
\]
Hence there is a map of spectral sequences
in hypercohomology. On the $E^2$ term, it is
 \[
\CD
\Ext^{j}_R(K_i,N) @>>> \clim^i\Ext^j(M_\lambda, N). 
\endCD
\]
The whole point is that the spectral sequence
on the left degenerates at $E^1$, since it
comes from a complex with zero differentials.
We have shown that
if $x\in\clim^i\Ext^{n+1-i}(M_\lambda, N)$
corresponds to an extension
\[
\CD
0 @>>> \y {\Sigma^{n+1-i}N}
@>>> F @>>> \y{K_{i-1}} @>>> 0
\endCD
\]
with $F\simeq\y Y$, then $x$ is the image of some $y$ from
the trivial spectral sequence. Therefore, all differentials
out of $x$ vanish.
\eprf

\exm{E2.11}
Let $k$ be an infinite field of cardinality $\aleph_t$.
Then by \cite{BaerLenzing82}, the polynomial ring
$R=k[x,y]$ is of pure global dimension $t+1$ ($\infty$ if $t$ is infinite). 
On the other hand, it is of global dimension $2$. Hence 
there do exist modules $N$  over $R=k[x,y]$,
satisfying the assumptions of the theorem when $t$ is at least 3.
\eexm

We can give a refinement of our results for when the ring $R$ is hereditary; 
recall that $R$ is {\em hereditary} if its global dimension is $\leq 1$. 
Examples of hereditary rings are commutative principal ideal domains,
and non-commutative polynomial rings.

\thm{th:hereditary}
Let $R$ be a hereditary ring.
Then
\begin{roenumerate}
\item {[BRM]} holds in $\ct$ if and only if the pure global dimension of $R$ 
is at most $1$; and
\item {[BRO]} holds in $\ct$ if and only if the pure global dimension of $R$ 
is at most $2$.
\end{roenumerate}
\ethm

\prf
(i) holds by Neeman's theorem~\ref{th:neeman}, combined with
the equality we prove in Proposition~\ref{P1.2}: for hereditary rings
\[
\pgldim R=\pgldim D(R).
\]

For (ii), note that Beligiannis' result (Theorem~\ref{T00.9})
tells us, that [BRO] holds if $\pgldim D(R)\leq 2$.
The converse comes from Proposition~\ref{P2.9} which says that if $N$ 
is an $R$--module of injective dimension $\leq1$ and $\PExt^3(M,N)\neq0$, 
then [BRO] fails for $\ct=D(R)$. Thus if $R$ is hereditary
but of pure global dimension $\geq3$, [BRO] must fail.
(Here we have used the easy fact that every hereditary ring is coherent.)
\eprf

Let $k$ be a field. In our counterexamples, we always consider
$k$-linear triangulated categories $\ct$.
When $\ct$ is $k$-linear, an additive functor $\{\ct^{c}\}^{op} \ri \ab$
always extends uniquely to a $k$-linear functor 
$\{\ct^{c}\}^{op}\ri\Mod{k}$, so we can restrict attention to such
$k$-linear functors.
The following lemma shows that our counterexamples must take
values in infinite-dimensional vector spaces.
The idea of the double dual used in the proof is due
to M.~Van den Bergh.

\lem{L2.12} 
Let $k$ be a field and 
\[
F:{\{\ct^c\}}^{op}\ri\mod{k}
\]
an exact functor which takes its values in the category
$\mod{k}$ of finite-dimensional vector spaces. Then
$F$ is of the form $\y{X}$ for some $X\in\ct$.
\elem

\prf Denote by $D$ the functor which takes a vector space 
to its dual. Then the functor
$G=D\circ F$ is exact and covariant. Let 
\[
\tilde{G}: \ct \ri \Mod{k}
\]
be the Kan extension of $G$ to $\ct$. Thus, for $Y\in\ct$,
we have
\[
\tilde{G}(Y) = \colim\, G C ,
\]
where the colimit is taken over the category of arrows
$C \ri Y$ from a compact $C$ to $Y$. A moment's thought
will convince the reader that $\tilde{G}$ is exact and
commutes with coproducts (cf. Prop.~2.3 of \cite{Krause98}). 
Hence $D\circ \tilde G$\/ is
exact and takes coproducts to products. By Brown's
theorem, it is representable: We have
\[
D\circ \tilde{G} = \ct\,(-, X)
\]
for some $X\in\ct$. We claim that $\y{X}= F$.
Indeed, the restriction of $D\circ \tilde{G}$\/
to $\ct^c$ is isomorphic to $D\circ D\circ F$,
and this functor is isomorphic to $F$ because
$FC$ is finite-dimensional for all $C\in\ct^c$.
\eprf

\bibliographystyle{amsplain}
\bibliography{stan}

\end{document}